\documentclass[11pt]{article} 
\usepackage[utf8]{inputenc}
\usepackage{caption,graphicx,amsmath,amsthm,amssymb, authblk,xcolor,enumitem,subfigure,float,verbatim,hyperref,pgfplots,bibentry,mathtools,soul,isomath, cancel}
\usepackage{mathrsfs}
\usepackage{indentfirst}
\usepackage{nccmath}
\graphicspath{ {images/} }
\usepackage{multicol, MnSymbol}
\usepackage{graphicx,amsmath,amssymb,amsthm,marvosym,tikz,fontenc}
\usetikzlibrary{backgrounds}
\usetikzlibrary{patterns, shapes.misc, positioning}
\usepackage{pgfplots}
\pgfplotsset{width=10cm,compat=1.9,tick scale binop=\times}
\usepackage{indentfirst}
\usepackage[a4paper,bindingoffset=0.2in,%
left=0.8in,right=1in,top=1.2in,bottom=1.2in,%
footskip=.25in]{geometry}
\usepackage{relsize}
\usepackage[english]{babel}
\allowdisplaybreaks
\theoremstyle{plain}
\newtheorem{theorem}{Theorem}[section]
\newtheorem{lemma}[theorem]{Lemma}
\newtheorem{proposition}[theorem]{Proposition}
\newtheorem{corollary}[theorem]{Corollary}
\newtheorem{definition}[theorem]{Definition}

\newtheorem{example}[theorem]{Example}

\date{}
\begin{document}
\title
{\bf{Eccentricity spectrum  of join of central graphs  and Eccentricity Wiener index of graphs}}
\author {\small Anjitha Ashokan
\footnote{anjithaashokan1996@gmail.com} \  and Chithra A V
\footnote{chithra@nitc.ac.in} \\ \small Department of Mathematics, National Institute of Technology, Calicut,\\
\small Kerala, India-673601}
\date{}
\maketitle
\begin{abstract} 

The eccentricity matrix of a simple connected graph is derived from its distance matrix by preserving the largest non-zero distance in each row and column, while the other entries are set to zero.
This article examines the $\epsilon$-spectrum, $\epsilon$-energy,  $\epsilon$-inertia and irreducibility of the central graph (respectively complement of the central graph)   of a triangle-free regular graph(respectively regular graph). Also look into the $\epsilon-$spectrum and the irreducibility of different central graph operations, such as central vertex join, central edge join, and central vertex-edge join. We also examine the $\epsilon-$ energy of some specific graphs. These findings allow us to construct new families of $\epsilon$-cospectral graphs and non $\epsilon$-cospectral $\epsilon-$equienergetic graphs. Additionally, we investigate certain upper and lower bounds for the eccentricity Wiener index of graphs. Also, provide an upper bound for the eccentricity energy of a self-centered graph.
 \\


\end{abstract}

\section{Introduction}
All graphs considered in this paper are undirected, finite, and simple. Let $G=(V(G),E(G))$ be a $(p,q)$ graph with vertex set $V(G)=\{v_{1},\ldots , v_{p}\}$ and edge set $E(G)=\{e_{1},\ldots ,e_{q}\}.$ 
The complement $\overline{G}$ of $G$ is a graph with the vertex set $V(G)$ and edges consisting of pairs of nonadjacent vertices in $G.$
The line graph $L(G)$ of a graph  $G$ is a graph whose vertex set is equal to the edge set of $G$, and two vertices in L(G) are adjacent if the corresponding edges in G share a common endpoint. The line graph of $L(G)$ is denoted as $L^{2}(G).$

The adjacency matrix $A(G)=(a_{i,j})$ of $G$ is a $p\times p$ matrix whose $(ij)^{th}$ entry $a_{i,j}$ is $1$ if the $i^{th}$ vertex  is adjacent with the $j^{th}$ vertex $(i.e, v_{i}\sim v_{j})$ of $G$ and $0$ otherwise. 
The eigenvalues of $A(G)$ are called $A-$eigenvalues of $G.$
The collection of all $A-$eigenvalues of $G$ together with their multiplicities are called $A-$spectra. The incidence matrix $R(G)=(r_{i,j})$ of $G$ is a $p\times q$ matrix whose $(ij)^{th}$ entry $r_{ij}$ is $1$ if the $i^{th}$ vertex of $G$ is incident with the $j^{th}$ edge of $G$ and $0$ otherwise. 
The distance matrix $D(G)=(d_{i,j})$ of a graph $G$ is a $p\times p$ matrix whose $(ij)^{th}$ entry $d_{i,j}$ represents the distance between the vertices $v_{i}$ and $v_{j}.$ The distance between the vertices $v_{i}$ and $v_{j}$ is denoted by $d(v_{i},v_{j}).$ The degree of a vertex, $deg(v_{i}),$ is the number of edges incident with vertex $v_{i}.$ 
The eccentricity of a vertex $v_{i}$ is defined as $e(v_{i})=\max \{d(v_{i},v_{j}):v_{j}\in V(G)\}.$ 
The vertex  $v_{k}$ is an eccentric vertex of $v_{i}$ if $d(v_{i}, v_{k}) = e(v_{i}).$
The minimum and maximum eccentricities of all vertices of $G$ are called radius($r(G)$) and diameter$(diam(G))$ respectively. A graph $G$ is self-centered if $r(G)$ and $diam(G)$ are the same. The total eccentricity of a graph is defined as $\varepsilon^{*}(G)=\sum_{i=1}^{p}\epsilon(v_{i}),$ and the eccentric connectivity index of $G$ is, $\zeta(G)=\sum_{i=1}^{p}deg(v_{i})e(v_{i}).$
\\
The eccentricity matrix, $\epsilon(G)$ of a graph $G$ is a $p\times p$ matrix its  entries are defined as   \begin{equation*}
    (\epsilon(G))_{ij}=\begin{cases}
            d(v_{i},v_{j}) & \text{ if } d(v_{i},v_{j})=\min\{e(v_{i}),e(v_{j})\},\\
            0 &\text{ otherwise.}
            \end{cases}
\end{equation*}
 A graph $G$ is $\epsilon-$regular if $\epsilon(i)=k,$ for all $i,$ where $\epsilon(i)=\sum_{j=1}^{n}\epsilon(G))_{ij}.$
Note that $\epsilon(G)$ is a real symmetric matrix, its eigenvalues are represented as $\epsilon-$eigenvalues of $G.$ If $\epsilon_{1}\geq \epsilon_{2} \geq , \ldots , \geq\epsilon_{k}$ are the distinct eigenvalues of $\epsilon(G)$  then $\epsilon-$spectrum of $\epsilon(G)$ is defined as 
\begin{align*}
spec_{\epsilon}(G)=\begin{pmatrix}
                    \epsilon_{1}& \epsilon_{2}&\cdots &\epsilon_{k}\\
                    m_{1}&m_{2}&\cdots &m_{k}
\end{pmatrix},               
\end{align*}
where $m_{i}$ denotes the algebraic multiplicity of the eigenvalue $\epsilon_{i}.$
The eccentric graph of   $G$, $G^{e}$  has the same vertex set as of  $G$. In $ G^{e}$, two vertices $ v_{i} $ and $ v_{j} $ are adjacent if and only if $ d(v_{i}, v_{j}) = \min\{e(v_{i}), e(v_{j})\}.$
Let $J$ and $I$ represent the all-one matrix and the identity matrix of appropriate order respectively and Q represents the equitable quotient matrix\cite{brouwer2011spectra}.

If a permutation matrix $P$ exists such that $$M=P^{T}\begin{pmatrix}
  N_{11} &N_{12} \\
  0 & N_{22}
\end{pmatrix}P,$$
where $N_{11}$ and $N_{22}$ are square block matrices, then $M$ is reducible. Otherwise, $M$ is irreducible. Unlike adjacency and distance matrices, the eccentricity matrix is not always irreducible. Some classes of graphs with irreducible or reducible eccentricity matrices are discussed in 
\cite{li2022inertia, MR4792389, MR4305906,MR3906706,MR4092627}

The notion of eccentricity matrix was introduced and studied in \cite{MR3906706}. The eccentricity matrix is also known as $D_{MAX}$ matrix in the literature\cite{MR3136762}. 
The $\epsilon-$spectral radius, $\rho_{\epsilon}(G)$ of $G$ is the largest value among the absolute values of $\epsilon-$eigenvalues and $\epsilon-$energy, $E_{\epsilon}(G)$ of $G$ is the sum of the absolute values of the $\epsilon-$eigenvalues.

Two non-isomorphic graphs with the same order and the same  $\epsilon-$spectrum are called $\epsilon-$cospectral graphs. Otherwise, they are non $\epsilon-$cospectral. Two graphs are $\epsilon-$equienergetic if they have the same $\epsilon-$ energy. Trivially all $\epsilon-$cospectral graphs are $\epsilon-$equienergetic. But the converse need not be true. Construction of non $\epsilon-$ cospectral $\epsilon-$equienergetic graphs is an interesting problem in spectral graph theory. In this article, we explore some families of non $\epsilon-$cospectral $\epsilon-$equienergetic graphs. The eccentricity Wiener index ($\epsilon-$ Wiener index) of a graph $G$ is defined as, $W_{\epsilon}(G)=\frac{1}{2}\sum_{ij} (\epsilon(G))_{ij}.$
Let $M$ be a symmetric matrix then the inertia of $M,$ $In(M)=(n_{+}(M), n_{-}(M), n_{0}(M))$ which represents the number of positive, negative, and zero eigenvalues of $M$, respectively. In literature, the inertia of eccentricity matrices of a few graphs are available, see 
\cite{MR4730406,li2022inertia, MR4110109, MR4446122, MR4305906}.
In this article, we discuss about the eccentricity inertia of some graphs.
\\
In graph theory, operations on graphs play an important role in constructing graphs from given graphs. See \cite{MR4792389, MR4305906}  for important studies on eccentricity spectra of graph operations. 
The central graph, $C[G],$ of a graph $G$ is a unary graph operation that is obtained from $G$ by adding new vertices, $I(G)$, to each edge of $G$ and joining all of its non-adjacent vertices. 
\\

In \cite{MR4657005, MR4780061, MR4150216}, the authors introduced various graph operations on central graphs and determined spectrum (adjacency, distance). Motivated by these findings, we examine the $\epsilon$-spectrum of central graphs, central vertex joins, central edge joins, and central vertex-edge joins of graphs. 
\\

In 
\cite{MR4574900, MR4727625}
there are some applications of the $\epsilon-$Wiener index as a lower bound for the eccentricity spectral radius and eccentricity spectral spread. Inspired by these studies, we establish the relationship between the eccentricity Wiener index and the total eccentricity of a graph.

This article is organized as follows: In Section 2, we gather some preliminary results. In Section 3, we obtain the $\epsilon-$spectrum, irreducibility, $\epsilon-$spectral radius, $\epsilon-$energy, and the inertia of the central graph of a graph.
Furthermore, the $\epsilon-$spectrum and $\epsilon-$energy of the complement of the central graph are estimated.
Moreover, the $\epsilon-$ spectrum, irreducibility, $\epsilon-$ Wiener index, and a lower bound for $\epsilon-$spectral radius are found for the central vertex join and central edge join of graphs, where the first graph is a regular triangle-free and the second graph is a regular graph.
In addition, the $\epsilon-$spectrum of the central vertex-edge join operation of a triangle-free regular graph with two regular graphs is computed.
Using these graph operations,  new families of $\epsilon-$ cospectral graphs and non-$\epsilon$-cospectral $\epsilon$-equienergetic graphs are constructed.
Additionally, we obtain the $\epsilon-$ energy of some particular class of graphs. 
Section 4 provides some lower and upper bounds for the $\epsilon-$ Wiener index of $G,$ also we provide a Nordhus-Guddum type upper bound for the $\epsilon-$ Wiener index. Moreover, we present an upper bound for the $\epsilon-$Wiener index for trees and for the $\epsilon-$ energy of the self-centered graph.

\section{Preliminaries}
This section contains some definitions and results that will be used in the next sections.
\begin{definition}\cite{MR4150216}
     Let $G_{i}$ be a $(p_{i},q_{i})$, $i=1,2$ graph.  The central vertex join of $G_{1}$ and $G_{2}$ is the graph $C[G_{1}]\dot{\vee} G_{2}$ is obtained by $C[G_{1}]$ and $G_{2}$ by joining each vertex of $G_{1}$ with every vertex of $G_{2}.$  
\end{definition}
\begin{definition}\cite{MR4150216}
    Let $G_{i}$ be a $(p_{i},q_{i})$, $i=1,2$ graph.  The central edge join of $G_{1}$ and $G_{2}$ is the graph $C[G_{1}]\veebar{G_{2}}$ is obtained by $C[G_{1}]$ and $G_{2}$ by joining each vertex corresponding to the edges of $G_{1}$ ($ I(G)$) with every vertex of $G_{2}.$  
\end{definition}

\begin{definition}\cite{MR4657005}
    Let $G_{1}, G_{2},G_{3}$ be any three graphs. The central vertex-edge join of $G_{1}$ with $G_{2}$ and $G_{3}$ is the graph $C[G_{1}] \vee (G_{2}^{V} \cup G_{3}^{E}),$  is obtained from $C[G_{1}],$ $G_{2}$ and $G_{3}$ by joining each vertex of $G_{1}$ with every vertex of $G_{2}$ and each vertex corresponding to the edges of $G_{1}$ with every vertex of $G_{3}.$
\end{definition}

\begin{lemma}\label{linegraphproB}\cite{MR0572262}
Let $G$ be an $r-$regular, $(p,q)$ graph with an adjacency matrix $A$ and an incidence matrix $R$. Let $L(G)$ be the line graph of $G$. Then $RR^{T}=A+rI$, $R^{T}R=B+2I$, where $B$ is the adjacency matrix of $L(G)$. Also, if $J$ is an all-one  matrix of appropriate order, then $JR=2J=R^{T}J$ and $JR^{T}=rJ=RJ.$ 
\end{lemma}
\begin{lemma}\label{spectrum of L(G)}\cite{MR0572262}
    Let $G$ be an $r-$regular, $(p,q)$ graph with $speC[G]=\{r,\lambda_{2},\ldots, \lambda_{p}\}$. Then 
                   \begin{equation*}
                       spec(L(G))=\begin{pmatrix}
                           2r-2 & \lambda_{2}+r-2 & \cdots&\lambda_{p}+r-2& -2\\
                           1 & 1 &\cdots& 1 & q-p
                       \end{pmatrix}.
                   \end{equation*}
 Also, $Z$ is an eigenvector corresponding to the eigenvalue $-2$ if and only if $RZ=0$.
\end{lemma}

\begin{theorem}\cite{MR4522999}\label{eccentric gh connected}
    Let $G$ be a  $(p,q)$ graph. Then the matrix $\epsilon(G)$ is irreducible if and only if $G^{e}$ is a connected graph.
\end{theorem}

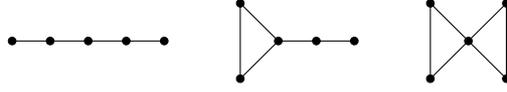
\begin{figure}[H]
 \centering
 \begin{tikzpicture}[scale=.5]
 \filldraw[fill=black](-6,0)circle(0.1cm); 
 \filldraw[fill=black](-5,0)circle(0.1cm);
 \filldraw[fill=black](-4,0)circle(0.1cm);
 \filldraw[fill=black](-3,0)circle(0.1cm);
 \filldraw[fill=black](-2,0)circle(0.1cm);
 \draw(-6,0)--(-5,0);
 \draw(-5,0)--(-4,0);
 \draw(-4,0)--(-3,0);
 \draw(-3,0)--(-2,0);
 \filldraw[fill=black](0,1)circle(0.1cm); 
 \filldraw[fill=black](0,-1)circle(0.1cm);
 \filldraw[fill=black](1,0)circle(0.1cm);
 \filldraw[fill=black](2,0)circle(0.1cm);
 \filldraw[fill=black](3,0)circle(0.1cm);
 \draw(0,1)--(1,0);
 \draw(0,1)--(0,-1);
 \draw(0,-1)--(1,0);
 \draw(1,0)--(2,0);
 \draw(2,0)--(3,0);
 \filldraw[fill=black](5,1)circle(0.1cm); 
 \filldraw[fill=black](5,-1)circle(0.1cm);
 \filldraw[fill=black](6,0)circle(0.1cm);
 \filldraw[fill=black](7,1)circle(0.1cm);
 \filldraw[fill=black](7,-1)circle(0.1cm);
 \draw(5,1)--(6,0);
 \draw(5,1)--(5,-1);
 \draw(5,-1)--(6,0);
 \draw(6,0)--(7,1);
 \draw(6,0)--(7,-1);
 \draw(7,1)--(7,-1);
  \end{tikzpicture}
  \caption{$F_{1}, F_{2} \text{ and }F_{3}$}
 \label{fig:3}
 \end{figure}

 \begin{lemma}\label{diameter of L(G)one}\cite{MR3100519}
     Let $G$ be a connected graph on $p\geq 3$ vertices. Then $diam(L(G))=1$ if and only if $G=K_{3}$ or $G=K_{1,p-1}.$
 \end{lemma}
\begin{lemma}\cite{MR3100519}\label{Fi theorem}
    Let $G$ be a connected graph, none of the three graphs $F_{1}$, $F_{2}$ and $F_{3}$ (see Figure \ref{fig:3}) are an induced subgraph of $G$ if and only if $diam(L(G))\leq 2.$
 \end{lemma}

\begin{lemma}\label{diameter2adjacency and eccentricity}\cite{MR4416682}
   Let $G$ be a $r-$regular graph with  $diam(G)=2$ and $e(v_{i})>1,$ for every vertex $v_{i}\in V(G)$, then 
 $$spec_{\epsilon}(G)=\begin{pmatrix}
    2(n-r-1)&-2(1+\lambda_{2})&\cdots& -2(1+\lambda_{n})\\
    1&1&\cdots &1
\end{pmatrix},$$
where $r=\lambda_{1}> \lambda_{2}\geq\ldots \geq \lambda_{n}$ are the $A-$eigenvalues of $G.$
    
\end{lemma}
\begin{lemma}\label{complement diameter 2}\cite{MR3167887}
If for any two adjacent vertices $u$ and $v$ of a graph  $G$, there is a third vertex $w$ that is not adjacent to either $u$ or $v$ (Property(\dagger)).
Then
\begin{enumerate}
    \item $\overline{G}$ is connected
    \item $diam(\overline{G})=2.$
\end{enumerate}
\end{lemma}
\begin{lemma}\label{spectra of c(g)}\cite{MR4150216}
    Let $G$ be a $r-$ regular graph on $n$ vertices and $r\geq \lambda_{2}\geq\ldots, \geq \lambda_{n}$ be the adjacency eigenvalues of $G$. Then $$spec(C[G])=\begin{pmatrix}
        0& \frac{(n-1-r)\pm \sqrt{(n-1-r)^{2}+8r}}{2} & \frac{-1-\lambda_{i}\pm \sqrt{(1+\lambda_{i}^{2}+4(\lambda_{i}+r}}{2}\\
        \frac{n(r-2)}{2} &1 &1
    \end{pmatrix}.$$
\end{lemma}
\begin{lemma}\cite{MR4574900}\label{eccentricity spectral radius bound}
     Let $G$ be a connected graph on $p$ vertices and $q$ edges. Then $\rho_{\epsilon}(\epsilon(G))\geq 2\frac{W_{\epsilon}}{p}$. Moreover, equality holds if and only if $G$ is $\epsilon-$ regular.  
\end{lemma}

\section{Eccentricity spectrum of central graph operations  of graphs}

This section explores the $\epsilon-$ spectral properties of the central graph, as well as the complement of the central graph of a graph. As an application, we construct new families of $\epsilon-$ cospectral and non-$\epsilon-$ cospectral $\epsilon-$equienergetic graphs.

\begin{theorem}
    Let $G$ be a triangle-free, $r-$regular$(r\geq 2)$, $(p,q)$ graph. Then the $\epsilon-$ spectrum of $\epsilon(C[G])$ consists of 
    
               \begin{enumerate}
                 \item    $3$ with multiplicity $q-p$
                 \item    $\frac{(3-(3r+\lambda_{i}(G)))\pm \sqrt{(5\lambda_{i}(G)+3(r-1))^{2}+16(\lambda_{i}(G)+r)}}{2}$, $i=2,\ldots,p$
                 \item     $\frac{(3q-4r+3)\pm \sqrt{(3q-8r+3)^{2}+16(p-2)(q-r)}}{2},$
                  \end{enumerate}
where $r=\lambda_{1}>\lambda_{2}(G)\geq \ldots \geq \lambda_{p}(G)$ are eigenvalues of $A(G)$.
\end{theorem}

\begin{proof}
By suitable labeling of vertices, the eccentricity matrix of  $C[G]$ is of the form, 
 $$\epsilon(C[G])=\begin{pmatrix}
                          2A(G)&2(J-R(G))\\
                          2(J-R(G))^{T}&3(J-I-B(G)
                          \end{pmatrix},$$ where $A(G),$
 $R(G)$ and $B(G)$ represents the adjacency matrix of $G,$ incidence matrix of $G,$ and adjacency matrix of line graph of $G$ respectively.

Let $U_{i}(G)$ (for $i=2,\ldots, p)$ and $V_{j}(G)$ (for $j=1,\ldots, q-p)$   be the eigenvectors of $A(G)$ and $B(G)$ corresponding to the eigenvalues $\lambda_{i}$ and $-2$ respectively. By using Lemma \ref{linegraphproB}  $R(G)^{T}U_{i}(G)$ is an eigenvector of $B(G)$ corresponding to the eigenvalue $\lambda_{i}(G)+r-2.$\\

Since $\epsilon(C[G])\begin{pmatrix}
                                                                                                 0\\
                                                                                                 V_{i}(G)
                                                                                              \end{pmatrix}=3\begin{pmatrix}
                                                                                                0 \\
                                                                                                V_{i}(G)
                                                                                                \end{pmatrix}$,
   $3$ (with multiplicity $q-p$) is an eigenvalue of $\epsilon(C[G]).$                                                                                              
 Now the equitable quotient matrix $Q$ of $\epsilon(C[G])$ is $\begin{pmatrix}
     2r&2(q-r)\\
     2(p-2)&3(q-2r+1)
 \end{pmatrix}$, and  the $2$ eigenvalues of $Q$ are  $\frac{(3q-4r+3)\pm \sqrt{(3q-8r+3)^{2}+16(p-2)(q-r)}}{2}$, 
 these will be two eigenvalues of  $\epsilon(C[G])$ as well.

 Now we will find the remaining eigenvalues of $\epsilon(C[G])$. Let $\zeta$ be an eigenvalue of $\epsilon(C[G])$, find a real number $s$ such that
 $\epsilon(C[G])\chi_{i}(G)=\zeta\chi_{i}(G),$ where  $\chi_{i}(G)=\begin{pmatrix}
     sU_{i}(G)\\
     R(G)^{T}U_{i}(G)
 \end{pmatrix}.$
 On solving the equation $\epsilon(C[G])\chi_{i}(G)=\zeta\chi_{i}(G)$, we get 
 \begin{align}
     2s\lambda_{i}(G)-2(\lambda_{i}+r)=\zeta s \notag\\
     -2s-3(\lambda_{i}(G)+r-1)=\zeta \label{s value}\\
     2s^{2}+s(\lambda_{i}(G)+3(\lambda_{i}(G)+r-1))-2(\lambda_{i}(G)+r)=0 \notag
 \end{align}
So, $s$ has two values, $s=\frac{-(5\lambda_{i}(G)+3(r-1))\pm \sqrt{(5\lambda_{i}(G)+3(r-1))^{2}+16(\lambda_{i}(G)+r)}}{4}$.\\
Hence from $(\ref{s value})$, $$\zeta=\frac{(3-(3r+\lambda_{i}(G))\pm \sqrt{(5\lambda_{i}(G)+3(r-1))^{2}+16(\lambda_{i}(G)+r)}}{2}, i=2, \ldots, p.$$
Thus, we obtained  the remaining $2(p-1)$ number of eigenvalues of $\epsilon(C[G]).$
\end{proof}

\begin{corollary}
Let $G$ be a triangle-free, $r-$regular ($r\geq 2)$, $(p,q)$ graph.
    Then the  $\epsilon-$ specctral radius of $C[G]$, 
    \begin{center} 
        $$\rho_{\epsilon}(C[G])=\frac{(3q-4r+3)+\sqrt{(3q-8r+3)^{2}+16(p-2)(q-r)}}{2}.$$
    \end{center}
\end{corollary}
\begin{corollary}
Let $G$ be a triangle-free, $r-$regular ($r\geq 2)$, $(p,q)$ graph.
    Then  $In(\epsilon(C[G]))$ is $(q,p,0).$
\end{corollary}
\begin{corollary}
    Let $G$ be a triangle-free, $r-$regular ($r\geq 2$), $(p,q)$ graph with $A-$eigenvalues  $\{r, \lambda_{2},\ldots, \lambda_{p}\}$. Then the eccentricity energy of $C[G]$ is, \begin{equation*}
        E_{\epsilon}(C[G])=3(q-p)+\sqrt{(3q-8r+3)^{2}+16(p-2)(q-r)}+\sum_{i=2}^{p}\sqrt{(5\lambda_{i}+3(r-1))^{2}+16(\lambda_{i}+r)}.
    \end{equation*}
\end{corollary}

The subsequent result indicates that $C[G]$ is a graph with property(\dagger). Estimating the eccentricity spectrum of $C[G]$ for a general graph is a complex process. However, using the following result, ascertaining the eccentricity spectrum of $\overline{C[G]}$ is quite easy.

\begin{lemma}\label{diam c(g) 2}
Let $G$ be a connected graph on $p$ $ (p\geq 3)$ vertices. 
If $u$ and $v$ are any two adjacent vertices in $C[G]$, then there will be a third vertex $w$ in $C[G]$ that is not adjacent to both $u$ and $v.$
\end{lemma}
\begin{proof}
Let $V(C[G])=V(G)\cup I(G)$, $u$ and $v$ be two adjacent vertices in $G$.\\
Case 1: Let $u,v \in V(G)$, 
and $P$ be a $u-v$ path in $G.$
If $P$ is of length $2$, select the internal vertex of $P$ as the required vertex. If the length of $P$ exceeds $2$, select the subdivision vertex between any two internal vertices of $P$ as the required vertex.\\
\noindent Case 2: 
Let $u \in V(G)$ and $v \in I(G)$. If $u$ is a pendant vertex of $G,$ then there exist a vertex in $I(G)$  as the desired vertex. Otherwise, there exists a vertex in $V(G)$ as the desired vertex.
\end{proof}

From Lemma \ref{diam c(g) 2}, \ref{diameter2adjacency and eccentricity} and \ref{spectra of c(g)} we obtain the following Corollaries. 
\begin{corollary}\label{inertia of complement of c(g)}
   Let $G$ be a $r-$regular graph on $p$ vertices. Then $$In(\overline{C[G]})=(p,p,\frac{p(r-2)}{2}).$$ 
\end{corollary}

\begin{corollary}
    Let $G$ be a $r-$regular graph on $p$ vertices. Then $$E_{\epsilon}(\overline{C[G]})= 2 E_{A}(C[G]).$$
\end{corollary}

\begin{theorem}
    Let $G$ be a $(p,q)$ graph and $V(Tr')$ be the collection of all vertices in $G$ that are not part of any triangle in $G$. Then \begin{enumerate}
        \item  $\epsilon(C[G])$ is irreducible if $G$ is triangle-free.
        \item  $\epsilon(C[G])$ is irreducible if $G$ is not triangle-free and $V(Tr')\neq \phi$.
    \end{enumerate}  
\end{theorem}
\begin{proof}
         Case 1: Let $G$ be a triangle-free graph. If  $G$ $\cong$ $K_{1}$, $K_{2},$ then clearly $(C[G])^{e}$ is connected. 
         If $G$ $\ncong$ $K_{1}$, $K_{2}$, 
        consider $C[G]$, with $V(C[G])=V(G)\cup I(G).$
        Let $u_{i}$ and $u_{j}'$ denotes the vertices in $V(G)$ and $I(G)$ respectively.\\
        In $C[G],$ \begin{align*}
          e(u_{i})=2, & \text{ for every } u_{i} \in V(G) \\
          2\leq e(u_{j}') \leq 3 & \text{ for every } u_{j}' \in I(G)
           \end{align*}
       If $u_{i}\sim u_{k}$ in $G$, then in $C[G]$, $d(u_{i},u_{k})=2=\min\{e(u_{i}), e(u_{k})\}.$ Therefore in $(C[G])^{e},$ $u_{i}$ is adjacent to $u_{k}.$ Since $G$ is a connected graph, the set of vertices in $V(G)$ is connected in $(C[G])^{e}.$ \\
       \noindent Now, let $u_{i}$ and $u_{k}$ be two adjacent vertices in $G,$ since $G$ is a triangle-free graph, there exists at least one vertex $u_{l}$ in $V(G)$ such that $u_{i}\sim u_{k}, u_{i}\sim u_{l}$ but $u_{k}\nsim u_{l}$. Then in $C[G]$, $u_{l}\sim u_{k}$.  If $u_{j}'$ is the subdivision vertex of the edge $(u_{i}, u_{k})$ then in $C[G]$, $d(u_{l}, u_{j}')=2=\min \{ e(u_{l}), e(u_{j}') \}.$ This implies every vertex $u_{j}'$ in $I(G)$  are adjacent to some $u_{l}$ in $V(G)$. Hence $(C[G])^{e}$ is connected. Therefore $\epsilon(C[G])$ is irreducible.\\
       \noindent Case 2: 
       Let $V(Tr)$ be the collection of vertices of $K_{3}$'s in $G$, $V(Tr')=V(G)\setminus V(Tr)$ , $V(I(Tr))$ be the collection of all newly added  vertices of $K_{3}$'s in $C[G]$, and $V(I(Tr'))=V(I(G))\setminus V(I(Tr)).$  
      If $V(Tr')\neq \phi$, then \begin{align*}
          e(v_{j})=2 & \text{ for every } v_{j} \in V(Tr')\\
          e(u_{i})=3 & \text{ for every } u_{i} \in V(Tr)\\
          3\leq e(u_{l}')\leq 4 & \text{ for every } u_{l} \in V(I(Tr))\\
          e(v_{k}')=3 & \text{ for every } v_{k}' \in V(I(Tr'))\\
      \end{align*}
Now, every vertices in $V(I(Tr))$ is an eccentric vertex of $v_{j},$ $v_{j}\in V(Tr')$. 
Therefore, $d(v_{j}, u_{l}')=2=\min\{e(v_{j}), e(u_{l}')\}$.
Hence, in $(C[G])^{e},$ $v_{j}$ is adjacent to every vertices in $V(I(Tr))$, for every vertices $v_{j}$ $\in$ $V(Tr')$.
Now clearly, for each $u_{i}$ in $V(Tr)$, there exist at least one $u_{l}'$ in $V(I(Tr))$ such that $d(u_{i}, u_{l}')=\min \{ e(u_{i}), e(u_{l}')\}$. Therefore, for every $u_{i} \in V(Tr)$ there exist atleat one $u_{l}'$ $\in$ $V(I(Tr))$ such that $u_{i}\sim u_{l}'$ in $(C[G])^{e}.$ Now the  eccentric vertices of $v_{k}'$ are in the set $V(I(Tr)).$ Thus, for each vertex $v_{k}'$ in $V(I(Tr'))$ there exist atleat one vertex $u_{l}'$ in $V(I(Tr))$ such that $d(v_{k}', u_{l}')=\min\{e(v_{k}'), e(u_{l}')\}$.
Therefore, for every $v_{k}' \in V(I(Tr'))$ there exist at least one $u_{l}'$ in $V(I(Tr))$ such that $v_{k}'\sim u_{l}'$ in 
$(C[G])^{e}$. Thus, $(C[G])^{e}$  is connected. By Theorem \ref{eccentric gh connected}  the matrix, $\epsilon(C[G])$ is irreducible.\\
\end{proof}
\begin{theorem}\label{centralvertexjoin}
For $i=1,2, $ let $G_{i}$ be a $r_{i}-$regular, $(p_{i},q_{i})$ graph, where $G_{1}$ is triangle-free. If  $\{r_{1}, \lambda_{2},\ldots \lambda_{p_{1}}\}$ and $\{r_{2}, \beta_{2},\ldots \beta_{p_{2}}\}$ are the $A-$ eigenvalues of $G_{1}$ and $G_{2}$ respectively.
Then the $\epsilon-$spectrum of $C[G_{1}]\dot{\vee} G_{2}$ consists of \begin{enumerate}
    \item $3$ with multiplicity $q_{1}-p_{1}$
    \item $-2(1+\beta_{j})$, $j=2,\ldots p_{2}$
    \item $-(2t_{1}+3(\lambda_{i}+r_{1}-1))$, $i=2,\ldots p_{1}$, $t_{1}=\frac{-(3(r_{1}-1)+5\lambda_{i})+\sqrt{(3(r_{1}-1)+5\lambda_{i})^{2}+16(\lambda_{i}+r_{1})}}{4}$
    \item $-(2t_{2}+3(\lambda_{i}+r_{1}-1))$, $i=2,\ldots p_{1}$, $t_{2}=\frac{-(3(r_{1}-1)+5\lambda_{i})-\sqrt{(3(r_{1}-1)+5\lambda_{i})^{2}+16(\lambda_{i}+r_{1})}}{4}$
    \item $3$ eigenvalues of the matrix $\begin{pmatrix}
        2r_{1}&2(q_{1}-r_{1})&0\\
        2(p_{1}-2)&3(q_{1}-2r_{1}+1)&2p_{2}\\
        0&2q_{1}&2(p_{2}-1-r_{2})
    \end{pmatrix}.$
\end{enumerate}
\end{theorem}
\begin{proof}
    By a proper labeling of vertices of $C[G_{1}]\dot{\vee} G_{2}$, $$\epsilon(C[G_{1}]\dot{\vee} G_{2})=\begin{pmatrix}
        2A(G_{1})& 2(J-R(G_{1}))&0\\
        2(J-R(G_{1})^{T}) & 3(J-I-B(G_{1}))&2J\\
        0&2J&2(J-I-A(G_{2}))
    \end{pmatrix},$$
    where $A(G_{i})$, $R(G_{1})$, and $B(G_{1})$ represents the adjacency matrix of $G_{i}(i=1,2)$, incidence matrix of $G_{1}$ and adjacency matrix of $L(G_{1})$ respectively.\\
Let $V_{j}(G_{1})(\text{ for } j=1,2\ldots q_{1}-p_{1})$ be an eigenvector of $B(G_{1})$ corresponding to the eigenvalue $-2$ (with multiplicity $q_{1}-p_{1}).$ Then we have $\epsilon(C[G_{1}]\dot{\vee} G_{2})\begin{pmatrix}
    0\\
    V_{j}(G_{1})\\
    0
\end{pmatrix}=3\begin{pmatrix}
    0\\
    V_{j}(G_{1})\\
    0
 \end{pmatrix}.$
Thus $3$ is an eigenvalue of $\epsilon(C[G_{1}]\dot{\vee} G_{2})$ with multiplicity $(q_{1}-p_{1}).$

Let $X_{k}(G_{2})$ be an eigenvector of $A(G_{2})$ corresponding to the eigenvalue $\beta_{k}$(for $k=2,\ldots, p_{2}$). Now $\epsilon(C[G_{1}]\dot{\vee} G_{2})\begin{pmatrix}
    0\\
    0\\
    X_{k}(G_{2})
\end{pmatrix}=-2(1+\beta_{j})\begin{pmatrix}
    0\\
    0\\
    X_{k}(G_{2})
 \end{pmatrix}.$
Thus $-2(1+\beta_{j})$(for $j=2,\ldots, p_{2}$) is an eigenvalue of $\epsilon(C[G_{1}]\dot{\vee} G_{2})$.

    Let $U_{l}(G_{1})$ be an eigenvector of $A(G_{1})$ corresponding the eigenvalue $\lambda_{l}$ (for $l=2,\ldots,p_{1})$.
    Using Lemma \ref{linegraphproB}, $R(G_{1})^{T}U_{l}(G_{1})$ is an eigenvector of $B(G_{1})$ corresponding to the eigenvalue $\lambda_{l}+r_{1}-2.$ Thus the vectors  $U_{l}(G_{1})$ and  $R(G_{1})^{T}U_{l}(G_{1})$ are orthogonal to the vector $J.$

     Now, consider the vector $\phi=\begin{pmatrix}
        tU_{l}(G_{1})\\
        R(G_{1})^{T}U_{l}(G_{1})\\
        0
    \end{pmatrix}.$
  Next,   we determine  under what condition $\phi$ is  an eigenvector of $\epsilon(C[G_{1}]\dot{\vee} G_{2})$. If $\mu$ is an eigenvalue of $\epsilon(C[G_{1}]\dot{\vee} G_{2})$ such that $\epsilon(C[G_{1}]\dot{\vee} G_{2})\phi= \mu \phi.$\\
  On solving we get, \begin{align*}
      2t\lambda_{i}-2(\lambda_{i}+r_{1})=&\mu t\\
      -(2t+3(\lambda_{i}+r_{1}-1))=&\mu
       \end{align*}
  From this we get,  
$$2t^{2}+t(3(\lambda_{i}+r_{1}-1)+2\lambda_{i})-2(\lambda_{i}+r_{1})=0$$
  
  Therefore, $\mu=-(2t+3(\lambda_{i}+r_{1}-1)),$ where $t=\frac{-(3(r_{1}-1)+5\lambda_{i})\pm \sqrt{(3(r_{1}-1)+5\lambda_{i})^{2}+16(\lambda_{i}+r_{1})}}{4}.$
  The remaining $3$ eigenvalues given by the equitable quotient matrix of $\epsilon(C[G_{1}]\dot{\vee} G_{2}),$ 
  $$Q=\begin{pmatrix}
        2r_{1}&2(q_{1}-r_{1})&0\\
        2(p_{1}-2)&3(q_{1}-2r_{1}+1)&2p_{2}\\
        0&2q_{1}&2(p_{2}-1-r_{2})
    \end{pmatrix}.$$
\end{proof}

\begin{corollary}
    Let $G_{i}$ be an $r_{i}-$regular, $(p_{i},q_{i})$ graph, $i=1,2.$ Then the  $\epsilon-$ Wiener index of  $\epsilon(C[G_{1}]\dot{\vee} G_{2})$ is given by 
    $$W_{\epsilon}(C[G_{1}]\dot{\vee} G_{2})=2q_{1}(p_{1}-1)+\frac{3}{2}q_{1}(q_{1}-2r_{1}+1)+p_{2}(2q_{1}+p_{2}-r_{2}-1).$$
\end{corollary}
\begin{corollary}
     Let $G_{i}$ be an $r_{i}-$regular, $(p_{i},q_{i})$ graph, $i=1,2.$  Then  
     $$\rho_{\epsilon}(C[G_{1}]\dot{\vee} G_{2})> \frac{4q_{1}(p_{1}-1)+3q_{1}(q_{1}-2r_{1}+1)+2p_{2}(2q_{1}+p_{2}-r_{2}-1)}{p_{1}+p_{2}+q_{1}}. $$ 
\end{corollary}
The following proposition can be obtained by using similar arguments as in Theorem 3.8.
\begin{proposition}\label{irreducibilityofvertex}
  Let $G_{1}$ and $G_{2}$ be any two graphs then $\epsilon(C[G_{1}]\dot{\vee} G_{2})$ is always irreducible. 
\end{proposition}

\begin{theorem}\label{central edge join}
For $i=1,2,$ let $G_{i}$ be a $r_{i}-$regular, $(p_{i},q_{i})$ graph where $G_{1}$ is triangle-free. If $\{r_{1}, \lambda_{2},\ldots \lambda_{p_{1}}\}$ and $\{r_{2}, \beta_{2},\ldots \beta_{p_{2}}\}$ are the $A-$ eigenvalues of $G_{1}$ and $G_{2}$ respectively, 
 then the $\epsilon-$ spectrum of $C[G_{1}]\veebar G_{2}$ consists of \begin{enumerate}
    \item $-2$ with multiplicity $q_{1}-p_{1}$
    \item $-2(1+\beta_{j})$, $j=2,\ldots p_{2}$
    \item $(\lambda_{i}-1)\pm \sqrt{(\lambda_{i}+1)^{2}+4(\lambda_{i}+r_{1})}$, $i=2, \ldots p_{1}$.
    \item $3$ eigenvalues of the matrix $\begin{pmatrix}
        2r_{1}&2(q_{1}-r_{1})&2p_{2}\\
        2(p_{1}-2)&2(q_{1}-1)&0\\
        2p_{1}&0&2(p_{2}-1-r_{2})
    \end{pmatrix}.$
\end{enumerate}
\end{theorem}
\begin{proof}
    By a proper labeling of vertices of $C[G_{1}]\veebar G_{2}$ we get \begin{align*}
      \epsilon(C[G_{1}]\veebar G_{2})=\begin{pmatrix}
       2A(G_{1}) & 2(J-R(G_{1}))& 2J\\
       2(J-R(G_{1}))^{T} &2(J-I)&0\\
       2J&0&2(J-I-A(G_{2}))
      \end{pmatrix}. 
    \end{align*}
    Let $Z_{l}(G_{1})(l=1,2\ldots q_{1}-p_{1})$ be an eigenvector of $B(G_{1})$(adjacency matrix of $L(G)$) corresponding to the eigenvalue $-2$(with multiplicity $q_{1}-p_{1}$). Then $ \epsilon(C[G_{1}]\veebar G_{2})\begin{pmatrix}
        0\\
        Z_{l}(G_{1})\\
        0
    \end{pmatrix}=-2\begin{pmatrix}
        0\\
        Z_{l}(G_{1})\\
        0
    \end{pmatrix}.$
    Thus $-2$ is an eigenvalue of $\epsilon(C[G_{1}]\veebar G_{2})$with multiplicity $q_{1}-p_{1}$.

    Let $X_{j}(G_{2})$ be an eigenvector of $A(G_{2})$ corresponding to the eigenvalue $\beta_{j}$ ($j=2,\ldots,p_{2}$). Then $ \epsilon(C[G_{1}]\veebar G_{2})\begin{pmatrix}
        0\\
        0\\
        X_{j}(G_{2})
    \end{pmatrix}=-2(1+\beta_{j})\begin{pmatrix}
        0\\
        0\\
        X_{j}(G_{2})
    \end{pmatrix}.$
    Thus $-2(1+\beta_{j})$ is an eigenvalue of $\epsilon(C[G_{1}]\veebar G_{2})$.

    Let $U_{k}(G_{1})$ be an eigenvector of $A(G_{1})$ corresponding to the eigenvalue $\lambda_{k}(k=2,\ldots,p_{1}).$ Then $R(G_{1})^{T}U_{k}(G_{1})$ is an eigenvector of $B(G_{1})$ corresponding to the eigenvalue $\lambda_{k}+r_{1}-2$. Consider the vector $\phi=\begin{pmatrix}
    tU_{k}(G_{1})\\
    R(G_{1})^{T}U_{k}(G_{1})\\
    0
    \end{pmatrix}$. Next, we will check the condition under which $\phi$ is an eigenvector of $\epsilon(C[G_{1}]\veebar G_{2}).$
If $\mu$ is an eigenvalue of $\epsilon(C[G_{1}]\veebar G_{2})$ such that $\epsilon(C[G_{1}]\veebar G_{2})\phi=\mu\phi$.
    On solving this we get $2$ values for $t$, $t=\frac{-(1+\lambda_{i})\pm \sqrt{(1+\lambda_{i})^{2}+4(r_{1}+\lambda_{i})}}{2}$.
    Thus, $\mu=(\lambda_{i}-1)\pm \sqrt{(\lambda_{i}+1)^{2}+4(\lambda_{i}+r_{1})}.$
    The remaining $3$ eigenvalues are given by the equitable quotient matrix of $\epsilon(C[G_{1}]\veebar G_{2}),$
    $$Q=\begin{pmatrix}
        2r_{1}&2(q_{1}-r_{1})&2p_{2}\\
        2(p_{1}-2)&2(q_{1}-1)&0\\
        2p_{1}&0&2(p_{2}-1-r_{2})
    \end{pmatrix}.$$
\end{proof}

\begin{corollary}
Let $G_{i}$ be an $r_{i}-$regular, $(p_{i},q_{i})$ graph for $i=1,2, $ where $G_{1}$ is triangle-free. Then the  $\epsilon-$ Wiener index of  $\epsilon(C[G_{1}]\veebar G_{2})$ is, 
    $$W_{\epsilon}(C[G_{1}]\veebar G_{2})=q_{1}(2p_{1}+q_{1}-3)+p_{2}(2p_{1}+p_{2}-1-r_{2}).$$
\end{corollary}
\begin{corollary}
     Let $G_{i}$ be an $r_{i}-$regular, $(p_{i},q_{i})$ graph for $i=1,2,$ where $G_{1}$ is triangle-free. Then  
     $$\rho_{\epsilon}(C[G_{1}]\veebar  G_{2})> \frac{2q_{1}(2p_{1}+q_{1}-3)+2p_{2}(2p_{1}+p_{2}-1-r_{2})}{p_{1}+p_{2}+q_{1}}. $$ 
\end{corollary}
\begin{theorem}
Let $G_{1}$ and $G_{2}$ be two graphs, then $\epsilon(C[G_{1}]\veebar G_{2})$ is always irreducible. 
\end{theorem}
\begin{proof}
    The proof follows from  a similar argument as in Theorem \ref{irreducibilityofvertex}.
\end{proof}

The following Corollaries are consequences of Theorems \ref{centralvertexjoin} and \ref{central edge join}.
\begin{corollary}
    Let $G$ be a regular, triangle-free graph and $S_{1},$  $S_{2}$ be two cospectral graphs. Then \begin{enumerate}
        \item $C[G]\dot{\vee}S_{1}$ and  $C[G]\dot{\vee}S_{2}$ are $\epsilon-$cospectral.
        \item $C[G]\veebar S_{1}$ and  $C[G]\veebar S_{2}$ are $\epsilon-$cospectral.
    \end{enumerate}
\end{corollary}
 \begin{corollary}
     Let $H_{1}$ and $H_{2}$ be two non-cospectral, $3-$regular, $(2t,3t)$  graph such that $t\geq 3$.
     Let $G_{1}=L^{2}(H_{1})$ and $G_{2}=L^{2}(H_{2})$. Let $G$ be a regular triangle-free graph. Then \begin{enumerate}
         \item $C[G]\dot{\vee}G_{1}$ and $C[G]\dot{\vee}G_{2}$ are non-$\epsilon$ cospectral $\epsilon-$ equienergetic.
         \item $C[G]\veebar G_{1}$ and $C[G]\veebar G_{2}$ are non-$\epsilon$ cospectral $\epsilon-$ equienergetic.
     \end{enumerate}
 \end{corollary}   
  \begin{theorem}
      For $i=1,2,3$, let $r_{i}=\lambda_{i1}\geq \lambda_{i2}\geq \ldots \lambda_{ip_{i}}$ be the $A-$eigenvalues of the $r_{i}-$regular  graph $G_{i}$, on $p_{i}$ vertices and $q_{i}$ edges, where $G_{1}$ is a triangle-free graph. Then the $\epsilon-$spectrum  of $C[G_{1}] \vee (G_{2}^{V} \cup G_{3}^{E})$ consists of, \begin{enumerate}
          \item $(-1+\lambda_{1j})\pm \sqrt{(1+\lambda_{1j})^{2}+4(\lambda_{1j}+r_{1})},j=2, 3, \ldots p_{1}.$
          \item $-2$ with multiplicity $q_{1}-p_{1}.$
          \item $0$ with multiplicity $p_{2}-p_{3}-2.$
          \item the eigenvalues of the matrix $\begin{pmatrix}
              2r_{1}&2(q_{1}-r_{1})&0&2p_{3}\\
              2(p_{1}-2)&2(q_{1}-1)&2p_{2}&0\\
              0&2q_{1}&0&3p_{3}\\
              2p_{1}&0&3p_{2}&0
               \end{pmatrix}.$
      \end{enumerate}
  \end{theorem}
  \begin{proof}
      By a proper labeling of vertices, $$\epsilon(C[G_{1}] \vee G_{2}^{V} \cup G_{3}^{E})=\begin{pmatrix}
          2A(G_{1})& 2(J-R(G_{1})) & 0& 2J\\
          2(J-(R(G_{1}))^{T})& 2(J-I)&2J & 0\\
          0&2J&0&J\\
          2J&0&3J&0
      \end{pmatrix}.$$
      The proof follows from a similar argument as in  Theorem \ref{centralvertexjoin}.
  \end{proof}
  
  \begin{corollary}
      Let $G_{1}$ be a triangle-free, regular graph and $G_{i}$, be $r-$regular, $(p,\frac{pr}{2})$ graph,, $i=2,3,4,5.$  Then $C[G_{1}] \vee (G_{2}^{V} \cup G_{3}^{E})$ and $C[G_{1}] \vee (G_{4}^{V} \cup G_{5}^{E})$ are $\epsilon-$ cospectral.
  \end{corollary}
\begin{theorem}
    Let $G_{i},$ $i=1,2,3$ be three graphs such that $G_{1}$ is triangle-free. Then $\epsilon(C[G_{1}] \vee G_{2}^{V} \cup G_{3}^{E})$ is irreducible.
\end{theorem}
\begin{proof}
    The proof follows from a similar argument as in Theorem \ref{irreducibilityofvertex}.
\end{proof}
  
\begin{theorem}\label{abc}
    Let $G$ be a $r-$regular  graph on $p$ $(p\geq 4)$ vertices such that none of the three graphs $F_{1},F_{2}$ and $F_{3}$ are induced subgraph of $G.$ If the smallest $A-$eigenvalue of $G$ is greater than or equal to $1-r,$ then $$E_{\epsilon}(L(G))=4p(r-1)+4(1-2r).$$
\end{theorem}
\begin{proof}
    Using Lemmas \ref{Fi theorem} and  \ref{diameter of L(G)one} and Theorem \ref{diameter2adjacency and eccentricity}, we get the desired result.
\end{proof}

\begin{corollary}
    Let $G_{1}$ and $G_{2}$ be a $r-$regular graphs on $p$ vertices, where $G_{1},G_{2}\neq K_{3}$, such that none of the three graphs $F_{1},$ $F_{2}$ and $F_{3}$ as an induced subgraph of $G_{i}, i=1, 2.$ Then $G_{1}$ and $G_{2}$ are cospectral if and only if $L(G_{1})$ and $L(G_{2})$ are $\epsilon-$ cospectral.
\end{corollary}
 \begin{corollary}
      let $G_{1}$ and $G_{2}$ be $r-$regular, non-complete , non-cospectral graphs on $p$ vertices such that none of the three graphs  $F_{1}$ $F_{2}$ and $F_{3}$  as an induced subgraph of $G_{i}, i=1, 2.$ If the smallest $A-$eigenvalue of $G_{i},i=1, 2$ greater than or equal to $1-r,$ then $L(G_{1})$ and $L(G_{2})$ are non $\epsilon-$ cospectral $\epsilon-$ equienergetic. 
 \end{corollary}
 The following Theorem is a consequence of 
   Lemmas \ref{spectrum of L(G)}, \ref{diameter2adjacency and eccentricity} and \ref{complement diameter 2}.
\begin{theorem}
    Let $G$ be a $r-$regular graph on   $p$ vertices and $L(G)$ be the line graph of $G.$
    If  the smallest $A-$eigenvalue of $G$ is greater than or equal to $2-r$ and $L(G)$ satisfy property (\dagger) then $$E_{\epsilon}(\overline{L(G)})=4p(r-2).$$
\end{theorem}

\begin{example}
 $CP(3)$ and $CP(4)$ are two  graphs with the smallest $A-$eigenvalue greater than or equal to $2-r$ and corresponding  line graphs satisfy property(\dagger).
\end{example}

\section{Some bounds for eccentricity Wiener index of graphs}
This section is devoted to determining some lower and upper bounds of the $\epsilon-$ Wiener index of graphs.
Let  $G$ is a $(p,q)$ graph with   $diameter$ $2$. Then  
\begin{align*}
W_{\epsilon}(G) = 
    \begin{cases}
               p(p-1)-2q+\frac{l}{2}(2p-l-1), &  \text{if  $G$ is not self centered}\\
              p(p-1)-2q, &  \text{if  $G$ is  self centered,}
            \end{cases} 
            \end{align*}
  where $l$ is the number of vertices $v_{i}$ such that $e(v_{i})=1.$  
 
\begin{proposition}\label{spectral radius of diameter 2 graph}
Let $G$ be a $(p,q)$ graph of diameter  $2.$ Then

\begin{equation*}
   \rho_{\epsilon}(G)\geq 
   \begin{cases}
			\frac{2(p(p-1)-2q+\frac{l}{2}(2p-l-1))}{p}, & \text{if  $G$ is not self centered}\\
           \frac{2(p(p-1)-2q)}{p}, & \text{if  $G$ is  self centered},
	\end{cases}
\end{equation*}
where $l$ is the number of vertices $v_{i}$ in $G$ such that $e(v_{i})=1.$  
\end{proposition}

\begin{theorem}\label{eccentricity wiener index for property graphs}
Let $G$ be a $(p,q)$ graph having property(\dagger). Then $W_{\epsilon}(\overline{G})=2q.$
    \end{theorem}
\begin{proof}
    The proof follows from Lemma \ref{complement diameter 2} and Proposition \ref{spectral radius of diameter 2 graph}.
\end{proof}
The following Corollary is a consequence of Theorem \ref{eccentricity wiener index for property graphs} and Lemma \ref{eccentricity spectral radius bound}.

\begin{corollary}
    Let $G$ be a  $(p,q)$ graph,  with girth $g(\geq 5)$. Then $$\rho_{\epsilon}(\overline{G})\geq \frac{4q}{p}.$$
\end{corollary}


Now, we consider graphs with  diameter greater than $2$. Determining the bounds for $W_{\epsilon}(G)$ based on the number of vertices and edges of $G$ is challenging. All upper bounds of the Wiener index of $G$ will also be upper bounds for $W_{\epsilon}(G)$. 
However, we obtain certain bounds for $W_{\epsilon}(G)$  that are tighter than the bounds of the Wiener index of $G$ by  using the total eccentricity of $G$ ($\varepsilon^{*}(G)$) and the eccentricity connectivity index($\zeta(G)$) of $G.$

 \begin{theorem}
  Let $G$ be a graph, then  $W_{\epsilon}(G)\geq \frac{\varepsilon^{*}(G)}{2}$. Equality holds if $G$ is an even cycle.
\end{theorem}

\begin{theorem}\label{helps for nordhus bound}
    Let $G$ be a $(p,q)$ graph with $e(v_{i})>1$ for every $v_{i} \in V(G)$. 
    Then \begin{align}
    W_{\epsilon}(G)\leq \frac{(p-1)\varepsilon^{*}(G)-\zeta(G)}{2}.\label{eque}
    \end{align}
    Equality holds if and only if $diam(G)=2.$
\end{theorem}
\begin{proof}
    From the definition of the eccentricity Wiener index we get $$2W_{\epsilon}(G)\leq \sum_{v_{i}}(p-1-deg(v_{i}))e(v_{i})=(p-1)\varepsilon^{*}(G)-\zeta(G).$$
    Moreover, we can easily see that equality holds if and only if $G$ is a graph having diameter $2.$
\end{proof}
Subsequently, using Theorem \ref{helps for nordhus bound} we provide a Nordhaus-Gaddum type upper bound for  $\epsilon-$ Wiener index of $G$.

\begin{theorem}
   Let $G$ be a graph $(p,q)$ with $e(v_{i})>1$, for every $v_{i}\in V(G)$. If $\overline{G}$ is also connected then $$W_{\epsilon}(G)+W_{\epsilon}(\overline{G})\leq \frac{1}{2}\Bigl((p-1)(\varepsilon^{*}(G)+\varepsilon^{*}(\overline{G}))-(\zeta(G)-\zeta(\overline{G}))\Bigr).$$
\end{theorem}
\begin{theorem}
    Let $T$ be a tree on $p$ vertices. Then 
    \begin{align*}\label{treewiemer}
        W_{\epsilon}(T)\leq \frac{k \varepsilon^{*}(T)+dk(p-k-1)-k(p-k)}{2},
    \end{align*} 
    where $k,d$ is the number of pendant vertices and the diameter of  $T$ respectively. Moreover, equality holds  if 
    $T \cong K_{1,p-1}.$ 
\end{theorem}
\begin{proof}
    Let $v_{1},\ldots, v_{k} $ be the pendant vertices of $T$ and $v_{k+1}, \ldots v_{p}$ be the non pendant vertices of $T.$
Then \begin{align*} 
\epsilon(v_{i}) & \leq  e(v_{i})(k-1)+(e(v_{i})-1)(p-k), \text{ for } i=1,2,\ldots, k.\\
\epsilon(v_{j}) & \leq  e(v_{j})k, \text{ for } j=k,k+1,\ldots, n.
  \end{align*}
Hence \begin{align*}
    2W_{\epsilon}(T) & \leq  k \varepsilon^{*}(T)+(p-k-1)(e(v_{1})+e(v_{2})\ldots e(v_{k}))-(p-k)k\\
    & \leq  k \varepsilon^{*}(T)+(p-k-1)kd-(p-k)k.
\end{align*}  
\end{proof}

Next, we provide an upper bound for the eccentricity energy of self-centered graphs in terms of its $\epsilon-$  Wiener index. 
\begin{theorem}
Let $G$ be a self-centered graph on $p$ vertices and diameter $d.$ Then
\begin{equation}\label{thm equation}
E_{\epsilon}(G)\leq \frac{2W_{\epsilon}(G)}{p} +\sqrt{2(p-1)W_{\epsilon}(G)\biggl(d-\frac{2W_{\epsilon}(G)}{p^{2}}\biggr)}.   
\end{equation}
 Equality holds in (\ref{thm equation}) if and only if $G$ is a $\epsilon-$regular graph with two  distinct $\epsilon-$eigenvalues $\frac{2W_{\epsilon}(G)}{p}$ and $-\sqrt{\frac{||\epsilon(G)||_{2}^{2}-\frac{4W_{\epsilon}(G)^{2}}{p^{2}}}{p-1}}$, or three distinct  $\epsilon-$eigenvalues $\frac{2W_{\epsilon}(G)}{p}$, $-\sqrt{\frac{||\epsilon(G)||_{2}^{2}-\frac{4W_{\epsilon}(G)^{2}}{p^{2}}}{p-1}}$ and $\sqrt{\frac{||\epsilon(G)||_{2}^{2}-\frac{4W_{\epsilon}(G)^{2}}{p^{2}}}{p-1}}$.
\end{theorem}
\begin{proof}
Let $\epsilon_{1}\geq \epsilon_{2}\geq \ldots \geq \epsilon_{p}$ be the $\epsilon-$eigenvalues of $G.$ Then by using Cauchy-Schwarz's inequality,  \begin{align}\label{cauchy}
    E_{\epsilon}(G)=\sum_{i=1}^{p}|\epsilon_{i}|&\leq \epsilon_{1}+\sqrt{(p-1)\sum_{i=2}^{p}|\epsilon_{i}|^{2}}\\
                                                &=\epsilon_{1}+\sqrt{(p-1)(||\epsilon(G)||_{2}^{2}-|\epsilon_{1}|^{2})}.\notag
\end{align}    
Now, consider the function $h(x)=x+\sqrt{(p-1)(||\epsilon(G)||_{2}^{2}-x^{2})},$ $h$ is strictly  decreasing in the interval $[\frac{||\epsilon(G)||_{2}}{\sqrt{p}},||\epsilon(G)||_{2}].$ 
Since $G$ is self-centered and by using Lemma \ref{eccentricity spectral radius bound} we have, $$\frac{||\epsilon(G)||_{2}}{\sqrt{p}}\leq \frac{2W_{\epsilon}(G)}{p}\leq \epsilon_{1}< ||\epsilon(G)||_{2}.$$ 
Therefore, \begin{equation}\label{monotone}
    \epsilon_{1}+\sqrt{(p-1)(||\epsilon(G)||_{2}^{2}-\epsilon_{1}^{2})}\leq \frac{2W_{\epsilon}(G)}{p}+\sqrt{(p-1)\biggl(||\epsilon(G)||_{2}^{2}-\frac{4W_{\epsilon}(G)^{2}}{p^{2}}\biggr)}.
\end{equation}
Hence, \begin{align*}
    E_{\epsilon}(G)   &\leq  \frac{2W_{\epsilon}(G)}{p}+\sqrt{(p-1)\biggl(||\epsilon(G)||_{2}^{2}-\frac{4W_{\epsilon}(G)^{2}}{p^{2}}\biggr)}\\
                      &=\frac{2W_{\epsilon}(G)}{p}+\sqrt{(p-1)\biggl( 2d W_{\epsilon}(G)-\frac{4W_{\epsilon}(G)^{2}}{p^{2}}\biggr)}.             
\end{align*}
Now we investigate the equality in (\ref{thm equation}). If the equality in (\ref{cauchy}) holds, then $|\epsilon_{2}|=|\epsilon_{3}|=\ldots=|\epsilon_{p}|=\sqrt{\frac{||\epsilon(G)||_{2}^{2}-\frac{4W_{\epsilon}(G)^{2}}{p^{2}}}{p-1}}.$ Equality in (\ref{monotone}) implies that   $\epsilon_{1}=\frac{2W_{\epsilon}(G)}{n}$.
 Then either $G$ is $\epsilon-$ regular with $2$ distinct $\epsilon-$ eigenvalues or $G$ is $\epsilon-$ regular with $3$ distinct $\epsilon-$ eigenvalues.  One can easily verify the converse.
This completes the proof.
\end{proof}

\section{Conclusion}
The graph operations such as central graph of a graph, central vertex join, central edge join and central vertex-edge join have been used to expand the collection of existing graphs for which  $\epsilon-$ spectrum, irreducibility, $\epsilon-$ inertia, and $\epsilon-$ energy are known. Also, the $\epsilon-$ energy of some classes of graphs are estimated.
These results allow us to construct new families of $\epsilon-$ cospectral graphs and non $\epsilon-$cospectral $\epsilon-$equienergetic graphs. Additionally, lower and upper bounds for the $\epsilon-$ Wiener index are determined by using the eccentricity connectivity index and the total eccentricity of a graph. Moreover,  an upper bound for the $\epsilon-$energy of self-centered graphs is provided.

\bibliographystyle{plain}
 \bibliography{reference}
\end{document}